\documentclass[11pt]{article}

\usepackage{amssymb}
\usepackage{amsmath}
\usepackage{amsfonts}
\usepackage[all]{xy}
\usepackage{graphicx}
\usepackage{epsfig}
\usepackage{makeidx}
\usepackage{color}
\input{epsf.tex}
\input xy
\xyoption{all}

\newtheorem{thm} {\noindent \bf{Theorem}}[section]

\newtheorem{remark}{\noindent \bf{Remark}}[section]

\newtheorem{prop}{\noindent \bf{Proposition}}[section]
\numberwithin{equation}{section}
\newenvironment{prooof}[1]{\noindent$\textbf{#1. }$}
{\hspace{\fill}$_\blacksquare$}

\begin{document}

\title{\sc Equilibrium balking strategies in the single server Markovian queue with
        catastrophes}
\author{Olga Boudali and Antonis Economou \\
        olboudali@math.uoa.gr and aeconom@math.uoa.gr\\
        University of Athens, Department of Mathematics\\
        Panepistemiopolis, Athens 15784, Greece}
\date{\today}
\maketitle

\noindent \textbf{Abstract:} We consider a Markovian queue subject
to Poisson generated catastrophes. Whenever a catastrophe occurs,
all customers are forced to abandon the system, the server is
rendered inoperative and an exponential repair time is set on. We
assume that the arriving customers decide whether to join the system
or balk, based on a natural reward-cost structure. We study the
balking behavior of the customers and derive the corresponding Nash
equilibrium strategies.

\vspace{0.5cm}


\vspace{1cm}

\section{Introduction}

Queues with removals of customers before being served are often
encountered in practice. One type of such a situation appears in
queueing systems with reneging, where customers are impatient and as
soon as their patience times expire they leave the system. Another
type of such a situation occurs in systems that are subject to
catastrophes/failures. Such events usually render the server(s)
inoperative and in addition force the customers to leave the
system. The crucial difference between these two types of
abandonments is that the customers decide whether to leave the
system or not according to their own desire in the case of reneging,
while they are forced to abandon the system in the case of catastrophes.\\

During the last decades, there is an emerging tendency to study
queueing systems from an economic viewpoint. More concretely, a
certain reward-cost structure is imposed on the system that reflects
the customers' desire for service and their unwillingness to wait.
Customers are allowed to make decisions about their actions in the
system, for example they may decide whether to join or balk, to wait
or abandon, to retry or not etc. The customers want to maximize
their benefit, taking into account that the other customers have the
same objective, and so the situation can be considered as a game
among the customers. In this type of studies, the main goal is to
find individual and social optimal strategies. The study of queueing
systems under a game-theoretic perspective was initiated by Naor
(1969) who studied the $M/M/1$ model with a linear reward-cost
structure. Naor (1969) assumed that an arriving customer observes
the number of customers and then makes his decision whether to join
or balk (observable case). His study was complemented by Edelson and
Hildebrand (1975) who considered the same queueing system but
assumed that the customers make their decisions without being
informed about the state of the system. Since then, there is a
growing number of papers that deal with the economic analysis of the
balking behavior of customers in variants of the $M/M/1$ queue, see
e.g. Burnetas and Economou (2007) ($M/M/1$ queue with setup times),
Economou and Kanta (2008a,b) ($M/M/1$ queue with compartmented
waiting space, $M/M/1$ queue with unreliable server), Guo and Zipkin
(2007) ($M/M/1$ queue with various levels of information and
non-linear reward-cost structure), Hassin and Haviv (1997) ($M/M/1$
queue with priorities), Hassin (2007) ($M/M/1$ queue with various
levels of information and uncertainty in the system parameters). The
monographs of Hassin and Haviv (2003) and Stidham (2009) summarize
the main approaches and several results in the broader area of the
economic analysis of queueing systems.\\

The study of the equilibrium customer behavior in queueing
systems with abandonments has received less attention. Hassin and
Haviv (1995) identified equilibrium customer strategies regarding
balking and reneging in the $M/M/1$ queue, where the reward for an
individual reduces to zero if its waiting time exceeds a certain
threshold time. Mandelbaum and Shimkin (2000) considered a quite
general model for abandonments from a queue, due to excessive wait,
assuming that waiting customers act rationally but without being
able to observe the queue length. More importantly they allowed
customers to be heterogeneous in their preferences and consequent
behavior. Other authors have also considered the equilibrium
behavior of customers in queueing systems with abandonments due to
reneging. Hassin and Haviv (2003), in Chapter 5, summarize the main
results for such models. However, to the best of our knowledge,
studies for the equilibrium behavior of customers in queueing
systems with abandonments/removals of customers due to catastrophic
events do not yet exist. It is the aim of the present
paper to study the equilibrium behavior of customers in the context
of a simple queueing model subject to catastrophes.\\

More specifically, in the present paper, we investigate the
equilibrium balking behavior of customers in a queue of $M/M/1$ type
with complete removals at Poisson generated catastrophe epochs. A
catastrophe renders the server inactive (due to either a failure or
preventive check/maintenance) and a repair time is set on. During
the repair time the system does not admit customers. When the repair
time is completed, the system behaves as an $M/M/1$ queue till the
next catastrophe and so on. We impose on the system a linear
reward-cost structure as the one in Naor (1969) and Edelson and
Hildebrand (1975). However, we make a modification, considering two
different types of reward: the first one is the usual reward
received by the customers that leave the system after service
completion, while the second is a compensation received by those
that are forced to abandon the system due to a catastrophe. In fact
the role of this compensation is to mitigate customers'
dissatisfaction. We again study the customers' behavior regarding
the dilemma whether to join or balk. We consider two cases with
respect to the level of information available to customers before
making their decisions. More specifically, at his arrival epoch, an
arbitrary customer may or may not be informed about the state of the
system (observable and unobservable cases correspondingly). In each
case, we characterize customer equilibrium strategies and we treat
the social optimization problem. We also explore the effect of the
information level on the equilibrium behavior of the customers
through numerical comparisons.\\

The paper is organized as follows. In Section 2 we describe the
dynamics of the model, the reward-cost structure and the decision
framework. In Section 3 we determine equilibrium threshold
strategies for the observable case, in which customers get informed
about the state of the system before making their decisions. In
Section 4 we study the unobservable case, deriving mixed equilibrium
balking strategies. Finally, in Section 5, we treat the social
optimization problem. Moreover, we present the results from several
numerical experiments that demonstrate the effect of the information
level on the behavior of the customers and on the
various performance measures of the system.\\

\section{Model description}

We consider a single-server queue with infinite waiting space, where
customers arrive according to a Poisson process at rate $\lambda$.
The service requirements of successive customers are independent and
identically distributed random variables with exponential
distribution with rate $\mu$. The server serves the customers one by
one. The system is subject to catastrophes/failures according to a
Poisson process at rate $\xi$. When a catastrophe occurs all
customers are forced to abandon the system prematurely, without
being served. The system is rendered inoperative and a repair
process is set on. The length of a repair time is exponentially
distributed at rate $\eta$. During a repair time, arrivals are not
accepted. We finally assume that interarrival times, service times,
intercatastrophe times and repair
times are mutually independent.\\

\begin{figure}[htp]
$$\def\labelstyle{\scriptstyle}
 \xymatrix {
(0,1)\ar@<1ex>[r]^{\lambda}
\ar@{<-}@<-1ex>[r]_{\mu}&(1,1)\ar@<1ex>[r]^{\lambda}
\ar@{<-}@<-1ex>[r]_{\mu}\ar@<1ex>[dl]^{\xi}
&(2,1)\ar@<1ex>[dll]^{\xi}\ar@<1ex>[r]^{\lambda}\ar@{<-}@<-1ex>[r]_{\mu}&\cdots\ar@<1ex>[r]^{\lambda}
&(n,1)\ar@<1ex>[dllll]^{\xi}\ar@<1ex>[r]^{\lambda}\ar@<1ex>[l]^{\mu}\ar@{<-}@<-1ex>[r]_{\mu}&\cdots\\
(0,0)\ar@<1ex>[u]^{\eta} \ar@{<-}@<-1ex>[u]_{\xi}}
$$\caption{Transition rate diagram of
$\{(Q(t),I(t))\}$.\label{f1}}
\end{figure}
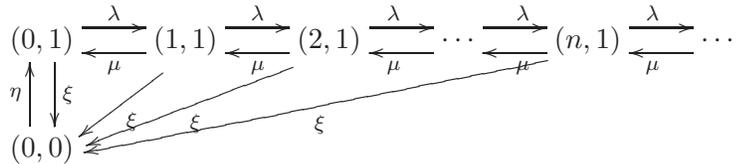

We represent the state of the system at time $t$ by a pair
$\left(Q(t),I(t)\right)$, where $Q(t)$ records the number of
customers at the system and $I(t)$ denotes the server state, with 1
describing a system in operation and 0 describing a system under
repair. Note that whenever $I(t)$ is zero, $Q(t)$ should be
necessarily zero too. Thus, the stochastic process
$\left\{\left(Q(t),I(t)\right):t\geq0\right\}$ is a continuous time
Markov chain with state space $S=\left\{\left(n, 1\right), n\geq
0\right\}\cup\left\{\left(0, 0\right)\right\}$ and its transition
rate diagram is shown in Figure \ref{f1}.\\

We are interested in the behavior of customers, when they have the
option to decide whether to join or balk. We model this decision
framework by assuming that each customer receives either a reward of
$R_{s}$ units for completing service or a compensation of $R_{f}$
units in case that he is forced to abandon the system due to a
failure. Moreover, a customer is charged a cost of $C$ units per
time unit that he remains in the system (either in queue or in the
service space). We also assume that customers are risk neutral and
wish to maximize their net benefit. Finally, their decisions are
assumed irrevocable, meaning that neither reneging of entering
customers nor retrials of balking customers are allowed.\\

Since all customers are assumed indistinguishable, we can consider
the situation as a symmetric game among them. Denote the common set
of strategies (set of available actions) and the payoff function by
$S$ and $F$ respectively. More concretely, let $F\left(s_{tagged},
s_{others}\right)$ be the payoff for a tagged customer who follows
strategy $s_{tagged}$, when all other customers follow $s_{others}$.
A strategy $s_{1}$ is said to dominate strategy $s_{2}$ if
$F\left(s_{1}, s\right)\geq F\left(s_{2}, s\right)$, for every $s
\in S$ and for at least one $s$ the inequality is strict. A strategy
$s_{*}$ is said to be weakly dominant if it dominates all other
strategies in $S$. A strategy $\tilde{s}$ is said to be a best
response against a strategy $s_{others}$, if $F\left(\tilde{s},
s_{others}\right)\geq F\left(s_{tagged},s_{others}\right)$, for
every $s_{tagged}\in S$. Finally, a strategy $s_{e}$ is said to be a
(symmetric) Nash equilibrium, if and only if it is a best response
against itself, i.e. $F\left(s_{e}, s_{e}\right)\geq F\left(s,
s_{e}\right)$, for every $s\in S$. The intuitive interpretation of a
Nash equilibrium is that it is a stable point of the game in the
sense that if all customers agree to follow it, then no one can
benefit by changing it. We remark that the notion of a dominant
strategy is stronger than the notion of an equilibrium. In fact,
every dominant strategy is an equilibrium, but the converse is not
true. Moreover, while equilibrium strategies
exist in most situations, a dominant strategy rarely does.\\

In the next sections we obtain customer equilibrium strategies for
joining/balking. We distinguish two cases with respect to the level
of information available to customers at their arrival instants,
before their decisions are made; the observable case (where
customers observe $Q\left(t\right)$) and the unobservable case.\\

\section{Equilibrium strategies - the observable case}

In this section we study the model, under the assumption that the
customers who find the server active observe the number of customers
in the system, before deciding whether to enter or balk. We prove
that a threshold type dominant strategy exists, which
constitutes the unique equilibrium balking strategy for the
customers in the system. We first give the expected net reward of a
customer that observes $n$ customers ahead of him and decides to
enter. We have the following.\\

\begin{prop}\label{customer-profit-mod1-obs}
Consider the observable model of the $M/M/1$ queue with catastrophes
causing complete removals of customers. The expected net benefit of
a customer that observes $n$ customers in the system upon arrival
and decides to enter is given by
\begin{equation}
S_{obs}(n)=R_{s}\left(\frac{\mu}{\mu+\xi}\right)^{n+1}+R_{f}
\left[1-\left(\frac{\mu}{\mu+\xi}\right)^{n+1}\right]
-
\frac{C}{\xi}\left[1-\left(\frac{\mu}{\mu+\xi}\right)^{n+1}\right],\;
n\geq 0.\label{S-obs1}
\end{equation}
\end{prop}

\begin{prooof}{Proof} Consider
a tagged customer that finds the system at state $(n,1)$ upon
arrival and decides to enter. This customer may leave the system
either due to its service completion or due to a catastrophe that
will force him to abandon prematurely the system. For his service
completion, he has to wait for a sum of $n+1$ independent
exponentially distributed times with parameter $\mu$ (note that
because of the memoryless property of the exponential distribution,
we can assume that the distribution of the remaining service time of
the customer in service is identical to the service time
distribution of the other customers). For the next catastrophe, he
has to wait for an exponentially distributed time with parameter
$\xi$. Therefore, the sojourn time of such a customer in the system
is given as $Z=\min (Y_n,X)$, where $Y_n$ follows a Gamma
distribution with parameters $n+1$, $\mu$ and $X$ is an
exponentially distributed random variable with rate $\xi$,
independent of $Y_n$. Moreover, the tagged customer will be served
with probability $\Pr [Y_n<X]$, while he will be forced to abandon
the system due to a catastrophe with the complementary probability
$\Pr [Y_n\geq X]$. Therefore his net benefit will be
\begin{equation}
S_{obs}(n)=R_{s}\Pr\left[Y_{n}<X\right]+R_{f}\Pr\left[Y_{n}\geq
X\right]-C E\left[Z\right].\label{S-obs1-initial}
\end{equation}
Note now that
\begin{equation}
\Pr\left[Y_{n}<X\right]=
\int_{0}^{\infty}e^{-\xi y}\frac{\mu^{n+1}}{n!} y^{n} e^{-\mu y}dy
=\left(\frac{\mu}{\mu+\xi}\right)^{n+1}\label{Pr[Y_n<X]}
\end{equation}
and
\begin{equation}
E\left[Z\right]=\int_{0}^{\infty}e^{-\xi z}\int_{z}^{\infty}\frac{\mu^{n+1}}{n!} u^{n} e^{-\mu u}du\ dz
=\frac{1}{\xi}\left[1-\left(\frac{\mu}{\mu+\xi}\right)^{n+1}\right].\label{E[Z]}
\end{equation}
Plugging \eqref{Pr[Y_n<X]} and \eqref{E[Z]} in \eqref{S-obs1-initial} yields \eqref{S-obs1}.
\end{prooof}\\

We now consider an arbitrary customer who observes upon arrival the
state of the system. Since arrivals are not permitted during the
repair time, if a customer observes the system at state $\left(0,
0\right)$, he is not allowed to enter and thus there is no decision.
So, we only consider the case where a customer observes the system
at a state $\left(n, 1\right)$. Such a customer strictly prefers to
enter if his net benefit is positive, is indifferent between joining
and balking if it is zero and strictly prefers to balk if it is
negative. In the sequel, we suppose for simplicity that customers
break ties in favor of entering. We have the following.

\begin{thm}\label{equilibrium-observable}
In the observable model of the $M/M/1$ queue with catastrophes
causing complete removals of customers, a unique
dominant pure strategy exists (which is also the unique equilibrium
strategy). There are three cases:
\begin{description}
\item[Case I:] $R_{f}<\frac{C}{\xi}-\frac{\mu R_{s}}{\xi}$.

Then the unique dominant strategy is always to balk.

\item[Case II:] $\frac{C}{\xi}-\frac{\mu R_{s}}{\xi}\leq
R_{f}<\frac{C}{\xi}$.

Then the unique dominant strategy is the threshold strategy `While
arriving at time $t$ and finding the system operative, observe
$Q(t)$; enter if $Q\left(t\right)\leq n_{e}$ and balk otherwise',
where $n_e$ is given by
\begin{equation}
n_e=\left\lfloor\frac{\ln K}{\ln S}-1\right\rfloor\label{n-e}
\end{equation}
with
\begin{equation}
K=\frac{\frac{C}{\xi}-R_{f}}{R_{s}-R_{f}+\frac{C}{\xi}},\;\;\;S=\frac{\mu}{\mu+\xi}
\label{K-S}
\end{equation}
and $\lfloor x \rfloor$ denotes the floor of $x$, i.e. the greatest
integer which is smaller than or equal to $x$.

\item[Case III:] $R_{f}\geq\frac{C}{\xi}$.

Then the unique dominant strategy is always to enter.
\end{description}
\end{thm}

\begin{prooof}{Proof} Consider a tagged customer that observes
the system upon arrival. If he finds the system at state $(n,1)$ and decides to enter,
then his expected net benefit is given by \eqref{S-obs1}. The
customer will prefer to enter if $S_{obs}(n)\geq 0$, which is
written easily as
\begin{equation}
\left(R_{s}-R_{f}+\frac{C}{\xi}\right)\cdot\left(\frac{\mu}{\mu+\xi}\right)^{n+1}\geq\frac{C}{\xi}-R_{f}.\label{1}
\end{equation}
Since $R_{s}-R_{f}+\frac{C}{\xi}>\frac{C}{\xi}-R_{f}$, we have the
following three cases.

Case A: $\frac{C}{\xi}-R_{f}>0\Leftrightarrow R_{f}<\frac{C}{\xi}.$

We can solve \eqref{1} with respect to $n$ and we obtain that the
tagged customer is willing to enter as long as he observes at most
$n_e$ customers in the system with $n_e$ given by \eqref{n-e}.
However, it is easy to see that $n_e$ given by \eqref{n-e} becomes
negative when $R_{f}<\frac{C}{\xi}-\frac{\mu R_{s}}{\xi}$.
Therefore, it is then optimal always to balk and we conclude with
Case I. On the other hand, when $\frac{C}{\xi}-\frac{\mu
R_{s}}{\xi}\leq R_{f}<\frac{C}{\xi}$, the threshold $n_e$ given
by \eqref{n-e} is non-negative and we conclude with Case II.

Case B:
$R_{s}-R_{f}+\frac{C}{\xi}>0\geq\frac{C}{\xi}-R_{f}\Leftrightarrow
\frac{C}{\xi}\leq R_{f}<\frac{C}{\xi}+R_{s}.$

In this case the inequality \eqref{1} is always true. Therefore the tagged
customer is always willing to enter.

Case C: $0\geq
R_{s}-R_{f}+\frac{C}{\xi}>\frac{C}{\xi}-R_{f}\Leftrightarrow
R_{f}\geq\frac{C}{\xi}+R_{s}.$

Solving inequality \eqref{1} with respect to $n$ shows that the
customer is willing to enter as long as as he observes at least
$n_e$ customers in the system with $n_e$ given by \eqref{n-e}.
However, $n_e$ is easily seen to be negative in this case, so it is
always preferable for the customer to enter. Therefore, Cases B
and C yield to Case III of the statement.

Note that the strategies prescribed above are dominant, since they
do not depend on what the other customers do, i.e. they are best
responses against any strategy of the others.
\end{prooof}

\begin{remark}
The equilibrium (dominant) strategies do not depend on the value of
the repair rate $\eta$. This happens because the customers make
decisions only whenever arrive at an operative system. On the
contrary, the social optimal strategies do depend on $\eta$, as we
will see in Section 5. Furthermore, in the limiting case where
$\xi\rightarrow0$, we can easily check, using L'Hospital rule, that
the threshold $n_e$ tends to the threshold derived by Naor (1969)
for the $M/M/1$ system.
\end{remark}

\section{Equilibrium strategies - the unobservable case}

We now turn our interest to the unobservable case, in which the
customers only know the values of the system parameters $\lambda$,
$\mu$, $\xi$ and $\eta$ and of the economic parameters $R_s$, $R_f$
and $C$, but do not observe the state of the system
upon arrival. Thus, there are only two pure strategies, `to join'
and `to balk' and a mixed strategy is specified by the joining
probability $q$ of an arriving customer that finds the server
operative. Our goal in this section is to identify
the equilibrium mixed balking strategies.\\

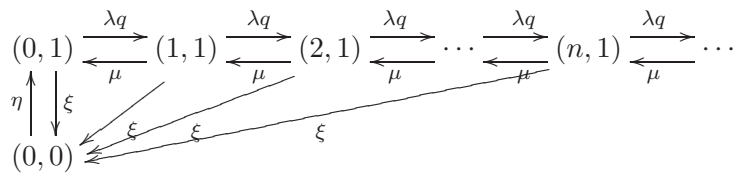
\begin{figure}[htp]
$$\def\labelstyle{\scriptstyle}
 \xymatrix {
(0,1)\ar@<1ex>[r]^{\lambda q}
\ar@{<-}@<-1ex>[r]_{\mu}&(1,1)\ar@<1ex>[r]^{\lambda q}
\ar@{<-}@<-1ex>[r]_{\mu}\ar@<1ex>[dl]^{\xi}
&(2,1)\ar@<1ex>[dll]^{\xi}\ar@<1ex>[r]^{\lambda
q}\ar@{<-}@<-1ex>[r]_{\mu}&\cdots\ar@<1ex>[r]^{\lambda q}
&(n,1)\ar@<1ex>[dllll]^{\xi}\ar@<1ex>[r]^{\lambda q}\ar@<1ex>[l]^{\mu}\ar@{<-}@<-1ex>[r]_{\mu}&\cdots\\
(0,0)\ar@<1ex>[u]^{\eta} \ar@{<-}@<-1ex>[u]_{\xi}}
$$\caption{Transition rate diagram of
$\{(Q(t),I(t))\}$ for a given mixed balking strategy $q$.\label{f3}}
\end{figure}

Suppose that the customers follow a mixed strategy with joining
probability $q$. Then, the system behaves as the original, but with
arrival rate $\lambda q$ instead of
$\lambda$. Its transition diagram is seen in Figure \ref{f3}.\\

We have the following.

\begin{prop}\label{stationary-unobservable}
Consider the unobservable model of the $M/M/1$ queue with
catastrophes causing complete removals of customers, in which the
customers that find the server operative join with probability $q$.
The stationary probabilities $p_{un}(k,i)$ of the system are given by
\begin{eqnarray}
p_{un}(0, 0)&=&\frac{\xi}{\xi+\eta},\label{7}\\
p_{un}(k,
1)&=&\frac{\eta\left(1-x_{2}(q)\right)x_{2}(q)^{k}}{\xi+\eta},\;\;\;
k\geq0,\label{7b}
\end{eqnarray}
where $x_2(q)$ is given by
\begin{equation}
x_{2}(q)=\frac{\left(\lambda q+\mu+\xi\right)-\sqrt{\left(\lambda
q+\mu+\xi\right)^{2}-4\lambda q\mu}}{2\mu}.\label{7c}
\end{equation}
The expected net benefit of a customer that enters with probability
$q'$ given that the system is found operative, when the others follow
a strategy $q$ is given by
\begin{equation}
S_{un}(q',q)=q' \left[ \left(R_{s}-R_{f}+\frac{C}{\xi}\right)
\frac{\mu\left(1-x_{2}(q)\right)}{\mu+\xi-\mu x_{2}(q)}
+R_{f}-\frac{C}{\xi} \right].\label{S-unobs}
\end{equation}
\end{prop}
\begin{prooof}{Proof} The balance equations for the stationary
distribution of the Markov chain  $\{(Q(t),I(t))\}$  are given as
follows:
\begin{eqnarray}
\eta p_{un}(0, 0)&=&\xi\sum_{k=0}^{\infty}p_{un}(k,
1),\label{5}\\
\left(\lambda q+\xi\right) p_{un}(0, 1)&=&\mu p_{un}(1,
1)+\eta p_{un}(0, 0),\label{4}\\
\left(\lambda q+\mu+\xi\right) p_{un}(k, 1)&=&\lambda q
p_{un}(k-1, 1)+ \mu p_{un}(k+1,
1),\;\;\;k\geq1.\label{6}
\end{eqnarray}
Equation \eqref{5} and the normalization equation $p_{un}(0,
0)+\sum_{k=0}^{\infty}p_{un}(k, 1)=1$ imply immediately \eqref{7}.
Equation \eqref{6} can be considered as a homogeneous linear
difference equation of order 2 with constant coefficients and
characteristic equation
\begin{equation}\label{11}
\left(\lambda q+\mu+\xi\right) x=\lambda q+\mu x^{2}\\
\end{equation}
that has two roots, $x_1(q)$ and $x_2(q)$, given by
\begin{equation}
x_{1,2}(q)=\frac{\left(\lambda
q+\mu+\xi\right)\pm\sqrt{\left(\lambda q+\mu+\xi\right)^{2}-4\lambda
q\mu}}{2\mu}.\label{x12-q}
\end{equation}
From the standard theory of homogeneous linear difference equations
(see e.g. Elaydi (1999) Section 2.3) we conclude that
$p_{un}(k,1)=c_1(q)x_1(q)^k+c_2(q)x_2(q)^k$, for $k\geq 0$, where
$c_1(q)$ and $c_2(q)$ are constants to be determined. We can easily
check that $x_1(q)>1$, hence $c_1(q)$ should be necessarily 0, for
$p_{un}(k,1),\; k\geq 0,$ are probabilities and so should remain
bounded. The constant $c_2(q)$ can be calculated using the
normalization equation and we deduce \eqref{7b}.

The expected net benefit of a customer that decides to enter when
the others follow the strategy $q$ can be computed by conditioning
on the state that he observes upon arrival. The probability that an
arriving customer finds $n$ customers in the system, given that he
finds the server operative (and so he can decide whether to enter or
not) is
\begin{equation}
p_{un}^{arr(\cdot,1)}(k,1)=\frac{\lambda
p_{un}(k,1)}{\sum_{i=0}^{\infty} \lambda
p_{un}(i,1)}=(1-x_2(q))x_2(q)^k,\;\;\; k\geq 0.
\end{equation}
Such a customer receives on the average $S_{obs}(k)$ units, given by
\eqref{S-obs1}. Therefore, the expected net benefit of a customer
that decides to enter given that he has found an operative system and
the others follow the strategy $q$ is
given by
\begin{eqnarray}
S_{un}(1,q)&=&\sum_{k=0}^{\infty} p_{un}^{arr(\cdot,1)}(k,1)
S_{obs}(k)\nonumber\\
&=& \sum_{k=0}^{\infty} (1-x_2(q))x_2(q)^k
\left\{R_{s}\left(\frac{\mu}{\mu+\xi}\right)^{k+1}+(R_{f}-\frac{C}{\xi})\left[1-\left(\frac{\mu}{\mu+\xi}\right)^{k+1}\right]
\right\}\nonumber\\
&=&\left(R_{s}-R_{f}+\frac{C}{\xi}\right)
\frac{\mu\left(1-x_{2}(q)\right)}{\mu+\xi-\mu x_{2}(q)}
+R_{f}-\frac{C}{\xi}.\label{S-un-1-q}
\end{eqnarray}
By the linearity of $S_{un}(q',q)$ with respect to the first
argument, we have that
$S_{un}(q',q)=(1-q')S_{un}(0,q)+q'S_{un}(1,q)$ and we obtain readily
\eqref{S-unobs}.
\end{prooof}\\

We can now proceed to determine the equilibrium balking strategies
of a customer in the unobservable case. We have the following.

\begin{thm}
In the unobservable model of the $M/M/1$ queue with catastrophes
causing complete removals of customers, a unique
equilibrium mixed strategy exists, with joining probability $q_e$ given by
\begin{equation}
q_e=
\begin{cases}
0&\text{if }R_{f}\leq\frac{C}{\xi}-\frac{\mu R_{s}}{\xi}\\
\frac{(C-\xi R_f+\xi R_s)(\mu R_s+\xi R_f-C)}{\lambda (C-\xi R_f)
R_s}& \text{if }\frac{C}{\xi}-\frac{\mu R_{s}}{\xi}< R_{f}<
\frac{C}{\xi}-\frac{\mu R_{s}\left(1-x_{2}\right)}{\xi}\\
1&\text{if }R_{f}\geq\frac{C}{\xi}-\frac{\mu
R_{s}\left(1-x_{2}\right)}{\xi},
\end{cases}\label{qe-definite}
\end{equation}
where $x_2=x_2(1)$ (using \eqref{7c} for $q=1$).
\end{thm}

\begin{prooof}{Proof} Suppose that customers who find the server
operative enter with probability $q$ and consider a tagged arriving
customer. Then, the tagged customer prefers to enter if
$S_{un}(1,q)>0$, he is indifferent between entering and balking if
$S_{un}(1,q)=0$ and he prefers to balk if $S_{un}(1,q)<0$. We
consider the equation $S_{un}(1,q)=0$ with $S_{un}(1,q)$ given by
\eqref{S-un-1-q} and we solve for $x_{2}(q)$. It may be easily checked
that the above equation has a unique solution given from
\begin{equation}
x_{2e}=\frac{\mu R_{s}+\xi R_{f}-C}{\mu R_{s}}\label{x-2e}
\end{equation}
and the corresponding $q_e$ is found by solving \eqref{11}, for
$x=x_{2e}$, with respect to $q$. This yields
\begin{equation}
q_e=\frac{x_{2e}[\mu (1-x_{2e})+\xi]}{\lambda
(1-x_{2e})}=\frac{(C-\xi R_f+\xi R_s)(\mu R_s+\xi R_f-C)}{\lambda
(C-\xi R_f) R_s}.\label{qe-general}
\end{equation}
We can now easily see that $q_e$ given by \eqref{qe-general} lies in
the interval $(0,1)$ if and only if $R_{f}\ \in
\left(\frac{C}{\xi}-\frac{\mu R_{s}}{\xi}, \frac{C}{\xi}-\frac{\mu
R_{s}\left(1-x_{2}\right)}{\xi}\right)$ and we obtain the second
branch of \eqref{qe-definite}. On the other hand, we can check that
$S_{un}(1,q)$ is positive for all $q\in [0,1]$, when
$R_{f}\geq\frac{C}{\xi}-\frac{\mu R_{s}\left(1-x_{2}\right)}{\xi}$
and, as a result, a customer's best response is 1 in this case.
Thus, `enter' is the unique equilibrium strategy and we obtain the
third branch of \eqref{qe-definite}.

Finally, for $q=0$, the system alternates between only two states,
$p_{un}(0, 0)$ and $p_{un}(0, 1)$, and the stationary probabilities are
$p_{un}(0, 0)=\frac{\xi}{\xi+\eta}$ and $p_{un}(0,
1)=\frac{\eta}{\xi+\eta}$. We can see that $S_{un}(1,0)$ is
non-positive if and only if $R_{f}\leq\frac{C}{\xi}-\frac{\mu
R_{s}}{\xi}$. Therefore, in this interval, $q=0$ is the unique
equilibrium strategy and we obtain the first branch of
\eqref{qe-definite}.
\end{prooof}

\begin{remark}
The equilibrium strategies do not depend on the value of the repair
rate $\eta$. This happens because the customers only make decisions
whenever they arrive at an operative system. However, unlike the
observable case, social optimal strategies do not depend on $\eta$
either, as we will see in Section 5. Furthermore, in the limiting
case where $\xi\rightarrow0$, we can easily check that the
equilibrium probability $q_e$ tends to the equilibrium probability
derived by Edelson and Hildebrand (1975) for the $M/M/1$ system.
\end{remark}

\section{Social optimal strategies - conclusions}

We are now studying the problem of social optimization. We treat
separately the observable and unobservable cases. We are interested
in determining the optimal values of the expected net social benefit
(per time unit) functions, $S_{obs}^{soc}(n)$ and
$S_{un}^{soc}(q)$, and the  corresponding arguments, $n_{soc}$ and
$q_{soc}$ respectively. First we consider the observable case and we
have the following Proposition \ref{Prop5.1}.

\begin{prop}\label{Prop5.1}
Consider the observable model of the $M/M/1$ queue with catastrophes
causing complete removals of customers. The expected net social
benefit per time unit, given that the customers follow a threshold
strategy with threshold $n$ (i.e. arriving customers that observe at
most $n$ customers in an operative system do enter, while the rest
balk without being served) is given by
\begin{eqnarray}
S_{obs}^{soc}(n)&=&\frac{\lambda(R_s-R_f)}{(\mu+\xi)^{n+1}} \left\{
\frac{\mu d_1 (n) \left[ (\mu+\xi)^{n+1}- \left( \mu x_1
\right)^{n+1} \right] }{\mu+\xi-\mu x_1} + \frac{\mu d_2 (n) \left[
(\mu+\xi)^{n+1}-
\left( \mu x_2 \right)^{n+1} \right] }{\mu+\xi-\mu x_2}\right\}\nonumber\\
& &+\lambda R_f \left( \frac{\eta}{\xi+\eta}-d_1 (n) x_1^{n+1} -d_2
(n)
x_2^{n+1} \right)\nonumber\\
& &-\frac{C\mu^2}{\xi^2} d_1 (n) x_1(1-x_2)^2 [1-(n+2)
x_1^{n+1}+(n+1)
x_1^{n+2}]\nonumber\\
& &-\frac{C\mu^2}{\xi^2} d_2 (n) x_2(1-x_1)^2 [1-(n+2)
x_2^{n+1}+(n+1) x_2^{n+2}],\; n\geq 0,\label{S-obs-soc}
\end{eqnarray}
where $x_1=x_1(1)$ and $x_2=x_2(1)$ (using \eqref{x12-q} for $q=1$)
and $d_1 (n)$, $d_2 (n)$ are given by
\begin{eqnarray}
d_1 (n)&=&\frac{-\eta \xi [(\mu+\xi)x_2-\lambda ]x_2^n}{(\xi+\eta)
\left\{ (\lambda+\xi-\mu x_2) [ (\mu+\xi) x_1 -\lambda ] x_1^n -
(\lambda+\xi-\mu x_1) [ (\mu+\xi) x_2 -\lambda ] x_2^n\right\}},\\
d_2 (n)&=&\frac{\eta \xi [(\mu+\xi)x_1-\lambda ]x_1^n}{(\xi+\eta)
\left\{ (\lambda+\xi-\mu x_2) [ (\mu+\xi) x_1 -\lambda ] x_1^n -
(\lambda+\xi-\mu x_1) [ (\mu+\xi) x_2 -\lambda ] x_2^n\right\}}.
\end{eqnarray}
\end{prop}
\begin{prooof}{Proof} The stationary distribution of the model under
a threshold strategy with threshold $n$ can be found along the same
lines with the proof of Proposition \ref{stationary-unobservable}
(i.e. by using the theory of linear difference equations with
constant coefficients). We then obtain that
\begin{eqnarray}
p_{obs}(0,0)&=&\frac{\xi}{\xi+\eta},\label{p00observable}\\
p_{obs}(k,1)&=&d_1 (n) x_1^k+d_2 (n) x_2^k,\;\;\;0\leq k \leq
n+1,\label{pk1observable}
\end{eqnarray}
with $x_1$, $x_2$, $d_1 (n)$ and $d_2 (n)$ as in the statement of
the Proposition. The expected net social benefit per time unit is
then found by
\begin{equation}
S_{obs}^{soc}(n)=\lambda P_{obs}^{ser} R_s +\lambda P_{obs}^{cat}
R_{f}-CE_{obs}[Q],\label{S-obs-soc-2}
\end{equation}
where $P_{obs}^{ser}$ and $P_{obs}^{cat}$ are the fraction of
customers that leave the system due to service and catastrophes
respectively and $E_{obs}[Q]$ is the mean number of customers in
system. Using \eqref{p00observable}, \eqref{pk1observable} and
\eqref{Pr[Y_n<X]}, we can compute $P_{obs}^{ser}$, $P_{obs}^{cat}$
and $E_{obs}[Q]$ as
\begin{eqnarray}
P_{obs}^{ser}&=&\sum_{k=0}^n  p_{obs}(k,1) \left( \frac{\mu}{\mu+\xi} \right)^{k+1}=\sum_{k=0}^n (d_1 (n) x_1^k +d_2 (n) x_2^k) \left( \frac{\mu}{\mu+\xi} \right)^{k+1},\label{Ps}\\
P_{obs}^{cat}&=&\sum_{k=0}^n  p_{obs}(k,1) \left[ 1 - \left( \frac{\mu}{\mu+\xi} \right)^{k+1} \right] \nonumber\\
&=&\sum_{k=0}^n (d_1 (n) x_1^k +d_2 (n) x_2^k) \left[ 1 - \left( \frac{\mu}{\mu+\xi} \right)^{k+1} \right], \label{Pf}\\
E_{obs}[Q]&=&\sum_{k=0}^{n+1} k p_{obs}(k,1)=\sum_{k=0}^{n+1} k (d_1
(n) x_1^k +d_2 (n) x_2^k).\label{EQ-observable}
\end{eqnarray}
Computing the relevant geometric sums in \eqref{Ps} - \eqref{EQ-observable} and substituting in
\eqref{S-obs-soc-2} yields \eqref{S-obs-soc}.
\end{prooof}\\

Unfortunately, the very involved form of \eqref{S-obs-soc} does not
allow the derivation of its maximum in closed analytic form.
However, it can be numerically evaluated quite easily. Thus, we turn
to numerical experiments below to derive some qualitative
conclusions for the behavior of the model.

\begin{figure}[hp]
\includegraphics[height=9cm,width=12cm,clip=true]{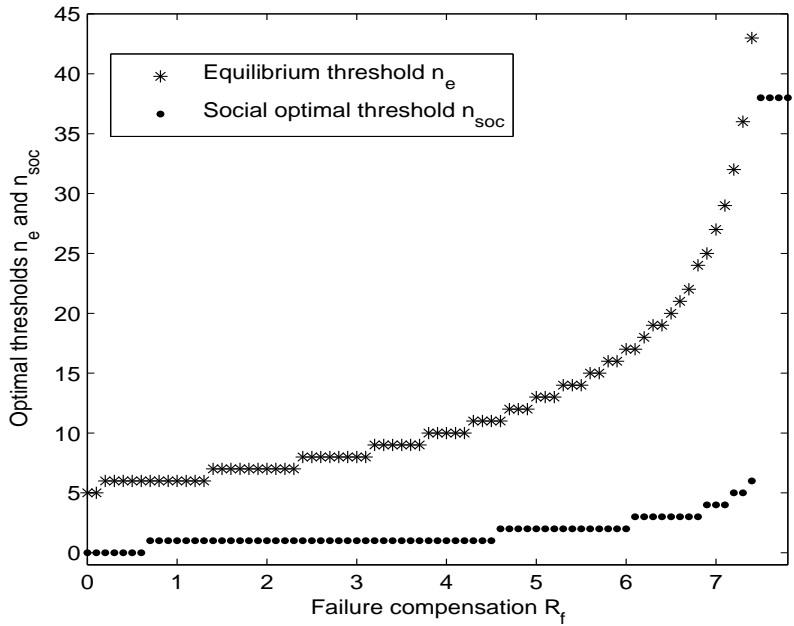}
\centering \caption{Equilibrium and social optimal joining
thresholds with respect to $R_f$ for the observable case with
$(\lambda,\mu,\xi,\eta)=(7,4,0.4,2)$ and $(R_s,C)=(7,3)$.}
\label{numerical-scenario1}\end{figure}

In Figure \ref{numerical-scenario1} we consider a model with
operation parameters $(\lambda,\mu,\xi,\eta)=(7,4,0.4,2)$ and
reward-cost parameters $(R_s,C)=(7,3)$ and we provide a graph of the
equilibrium and social optimal thresholds for the observable case as
functions of the failure compensation $R_f$. We observe that $n_e$
becomes infinity for large values of $R_f$, while $n_{soc}$
stabilizes to a certain value for large values of $R_f$. Moreover,
we observe that $n_{soc}\leq n_e$ for all values of $R_f$. These
qualitative facts seem to be valid in general, as it has been
verified from a large number of similar numerical experiments
for other values of the parameters.\\

We now turn to the unobservable case and we obtain the following
Proposition \ref{Prop5.2}.

\begin{prop}\label{Prop5.2}
Consider the unobservable model of the $M/M/1$ queue with
catastrophes causing complete removals of customers. The expected
net social benefit per time unit, given that the customers follow a
mixed strategy with joining probability $q$ (i.e. arriving customers
that find an operative system enter with probability $q$, while the
rest balk without being served) is given by
\begin{equation}
S_{un}^{soc}(q)=\frac{\eta x_2(q) [ \mu R_s (1-x_2(q))+\xi R_f -C]}{(\xi+\eta)(1-x_2(q))},\label{S-unobs-soc}
\end{equation}
with $x_2(q)$ given by \eqref{7c}.
\end{prop}
\begin{prooof}{Proof} The expected net social benefit per time unit is given by
\begin{equation}
S_{un}^{soc}(q)=\lambda P_{un}^{ser} R_s +\lambda P_{un}^{cat}
R_{f}-CE_{un}[Q],\label{S-unobs-soc-2}
\end{equation}
where $P_{un}^{ser}$ and $P_{un}^{cat}$ are the fractions of
customers that join but leave the system due to service and
catastrophes respectively and $E_{un}[Q]$ is the mean number of
customers in system. Using \eqref{7}, \eqref{7b} and
\eqref{Pr[Y_n<X]}, we can compute $P_{un}^{ser}$, $P_{un}^{cat}$ and
$E_{un}[Q]$ as
\begin{eqnarray}
P_{un}^{ser}&=&\sum_{k=0}^{\infty} p_{un}(k,1) q \left( \frac{\mu}{\mu+\xi} \right)^{k+1} =\sum_{k=0}^{\infty} \frac{\eta\left(1-x_{2}(q)\right)x_{2}(q)^{k}}{\xi+\eta} q \left( \frac{\mu}{\mu+\xi} \right)^{k+1} ,\label{Ps-unobs}\\
P_{un}^{cat}&=&\sum_{k=0}^{\infty} p_{un}(k,1) q  \left[ 1 - \left( \frac{\mu}{\mu+\xi} \right)^{k+1} \right] \nonumber\\
&=&\sum_{k=0}^{\infty} \frac{\eta\left(1-x_{2}(q)\right)x_{2}(q)^{k}}{\xi+\eta} q \left[ 1 - \left( \frac{\mu}{\mu+\xi} \right)^{k+1} \right], \label{Pf-unobs}\\
E_{un}[Q]&=&\sum_{k=0}^{\infty} k p_{un}(k,1)=\sum_{k=0}^{\infty} k
\frac{\eta\left(1-x_{2}(q)\right)x_{2}(q)^{k}}{\xi+\eta}
.\label{EQ-unobservable}
\end{eqnarray}
Computing the geometric sums in \eqref{Ps-unobs} -
\eqref{EQ-unobservable} and substituting in \eqref{S-unobs-soc-2}
yields \eqref{S-unobs-soc}.
\end{prooof}\\

Unlike the observable case, it is possible here to obtain the
optimal joining probability $q_{soc}$ in closed form and we can also
compare it with the corresponding equilibrium probability $q_{e}$.
The results are summarized in the following Theorem
\ref{Theorem5.1}.

\begin{thm}\label{Theorem5.1}
In the unobservable model of the $M/M/1$ queue with catastrophes
causing complete removals of customers, a unique social optimal
strategy exists, with joining probability $q_{soc}$ given by
\begin{equation}
q_{soc}=
\begin{cases}
0&\text{if }R_{f}\leq\frac{C}{\xi}-\frac{\mu R_{s}}{\xi}\\
\frac{\sqrt{D}\left(\mu R_{s}-\sqrt{D}\right)\left(\xi R_{s}+\sqrt{D}\right)}
{\lambda D R_{s}}& \text{if }\frac{C}{\xi}-\frac{\mu R_{s}}{\xi}< R_{f}<
\frac{C}{\xi}-\frac{\mu R_{s}\left(1-x_{2}\right)^{2}}{\xi}\\
1&\text{if }R_{f}\geq\frac{C}{\xi}-\frac{\mu R_{s}\left(1-x_{2}\right)^{2}}{\xi},
\end{cases}\label{qsoc-definite}
\end{equation}
where
\begin{equation}
 D=\mu R_{s}\left(C-\xi R_{f}\right)\label{D}
\end{equation}
and $x_2=x_2(1)$ (using \eqref{7c} for $q=1$). Moreover, the social
optimal joining probability is always smaller than the individual
one, i.e.
\begin{equation}
q_{soc}\leq q_{e}.\label{q-e-q-soc-inequality}
\end{equation}
\end{thm}

\begin{prooof}{Proof}
We observe that the function $S_{un}^{soc}(q)$ given by
\eqref{S-unobs-soc} can be written as the composition of $f(x)$ and
$x_2(q)$ (i.e. $S_{un}^{soc}(q)=f(x_2(q))$), with
\begin{equation}
f(x)=\frac{\eta x [ \mu R_s (1-x)+\xi R_f
-C]}{(\xi+\eta)(1-x)}\label{f-function}
\end{equation}
and $x_2(q)$ given by \eqref{7c}. Note also that the function
$x_2(q)$ is strictly increasing for $q\in [0,1]$ and it takes values
in $[0,x_2]$.

To proceed, we solve the equation
\begin{equation}
S_{un}^{soc\; '}(q)=f'(x_2(q))x_2'(q)=0,\label{S-un-soc-deriv=0}
\end{equation}
for $q\in [0,1]$. However, $x_2'(q)\ne 0$ for $q\in [0,1]$ and
therefore \eqref{S-un-soc-deriv=0} is reduced to $f'(x_2(q))=0$. So we
have to solve $f'(x)=0$,
which is also written after some straightforward algebra in the form
\begin{equation}
\mu R_s x^2 -2 \mu R_s x+ (\mu R_s +\xi R_f -C)=0.\label{f-deriv=0}
\end{equation}
The discriminant of the quadratic polynomial in \eqref{f-deriv=0} is
non-positive if and only if $R_f\geq \frac{C}{\xi}$. In that case,
we conclude that $f(x)$ is increasing and consequently
$S_{un}^{soc}(q)$ is also increasing. In summary we have:
\begin{itemize}
\item Case I: $R_f\geq \frac{C}{\xi}$. The social optimal joining
probability is $q_{soc}=1$.
\end{itemize}

In case where $R_f< \frac{C}{\xi}$, the equation $f'(x)=0$
(equivalently \eqref{f-deriv=0}) has two distinct roots $x_2^-$ and
$x_2^+$ given by
\begin{equation}
x_2^{-}=1- \frac{\sqrt{D}}{\mu R_{s}},\;\;\; x_2^{+}=1+ \frac{\sqrt{D}}{\mu R_{s}},\label{x-2+-}
\end{equation}
with $D$ given by \eqref{D}. Therefore, the quadratic polynomial in
\eqref{f-deriv=0} is positive for $x<x_2^-$ or $x>x_2^+$ and
negative for $x_2^-<x<x_2^+$. Due to the one to one correspondence
between $x_{2}(q)$ and $q$ and the fact that $x_2(q)\in [0,x_2]$ for
$q\in [0,1]$, we have to consider several cases regarding the
relative order of $x_2^-$, $x_2^+$ and $0$, $x_2$. However, we have
that $x_2<1<x_2^+$ and therefore there are only 3 subcases.
\begin{itemize}
\item Case II-a: $R_f< \frac{C}{\xi}$ and $x_2^-\leq 0$. Then, we have necessarily
$x_2^-\leq 0<x_2<x_2^+$. Then $f(x)$ is decreasing in $[0,x_2]$ and
consequently $S_{un}^{soc}(q)$ is decreasing in $[0,1]$. The social
optimal joining probability is $q_{soc}=0$.

\item Case II-b: $R_f< \frac{C}{\xi}$ and $0< x_2^- < x_2$.
Then, we have that $f'(x)$ is positive in $(0,x_2^-)$ and negative
in $(x_2^-,x_2)$ and therefore we conclude that
the maximum of $S_{un}^{soc}(q)$ is attained for $q$ such that
$x_2(q)=x_2^-$. The social optimal joining probability is found by
substituting $x_2^-$ for $x$ in \eqref{11} and solving for $q$. We obtain
\begin{equation}
q_{soc}=\frac{x_2^-[\mu (1-x_2^-) +\xi]}{\lambda (1-x_2^-)}
\end{equation}
and using \eqref{x-2+-} we deduce that
\begin{equation}
q_{soc}=\frac{\sqrt{D}\left(\mu R_{s}-\sqrt{D}\right)\left(\xi R_{s}+\sqrt{D}\right)}
{\lambda D R_{s}}.\label{qsoc-general}
\end{equation}

\item Case II-c: $R_f< \frac{C}{\xi}$ and $x_2^- \geq x_2$. Then we
have that $f(x)$ is increasing in $[0,x_2]$ and consequently $S_{un}^{soc}(q)$ is increasing in $[0,1]$. The social optimal joining probability is $q_{soc}=1$.
\end{itemize}
Using \eqref{x-2+-} and taking into account the common condition
$R_f< \frac{C}{\xi}$ for Cases IIa-c, we can easily see that the
conditions $x_2^-\leq 0$, $0< x_2^- < x_2$ and $x_2^- \geq x_2$ can
be written respectively  as $R_f \leq \frac{C}{\xi}-\frac{\mu
R_s}{\xi}$, $\frac{C}{\xi}-\frac{\mu
R_s}{\xi}<R_f<\frac{C}{\xi}-\frac{\mu R_s (1-x_2)^2}{\xi}$ and
$R_f\geq\frac{C}{\xi}-\frac{\mu R_s (1-x_2)^2}{\xi}$. By combining
Cases I-II we obtain immediately \eqref{qsoc-definite}.

Regarding the order between $q_{soc}$ and $q_{e}$, formulas
\eqref{qe-definite} and \eqref{qsoc-definite} show that
$q_e=q_{soc}=0$ for $R_f \leq \frac{C}{\xi}-\frac{\mu R_s}{\xi}$,
whereas $q_e=q_{soc}=1$ for $R_f\geq\frac{C}{\xi}-\frac{\mu R_s
(1-x_2)^2}{\xi}$. Moreover, $q_e=1$ and $q_{soc}\in (0,1)$ when
$\frac{C}{\xi}-\frac{\mu R_s (1-x_2)}{\xi}\leq
R_f<\frac{C}{\xi}-\frac{\mu R_s (1-x_2)^2}{\xi}$. Thus we have to
check the validity of the inequality $q_{soc}\leq q_e$, only for
$R_f\in (\frac{C}{\xi}-\frac{\mu R_s}{\xi},\frac{C}{\xi}-\frac{\mu
R_s (1-x_2)}{\xi})$, i.e. in the interval where both $q_e$ and
$q_{soc}$ are strictly between 0 and 1. Using \eqref{qe-general} and
\eqref{qsoc-general}, we can easily see after some
algebra that the inequality $q_{soc}\leq q_{e}$ is reduced to $\mu^2 \xi R_s^3+D\sqrt{D}\geq 0$ which clearly holds. Thus, the inequality is also valid in this case. \end{prooof}\\

The inequality $n_{soc}\leq n_e$ that the numerical experiments
suggest in combination with the inequality $q_{soc}\leq q_e$ that
has been analytically proved shows that we have the usual situation
also encountered in the pioneering papers of Naor (1969) and Edelson
and Hildebrand (1975): Individual optimization leads to longer
queues than it is socially desirable (i.e. in equilibrium the
customers make excessive use of the system). Indeed, a customer that
decides to join the system imposes negative externalities on future
arrivals.

\begin{figure}[hp]
\includegraphics[height=9cm,width=12cm,clip=true]{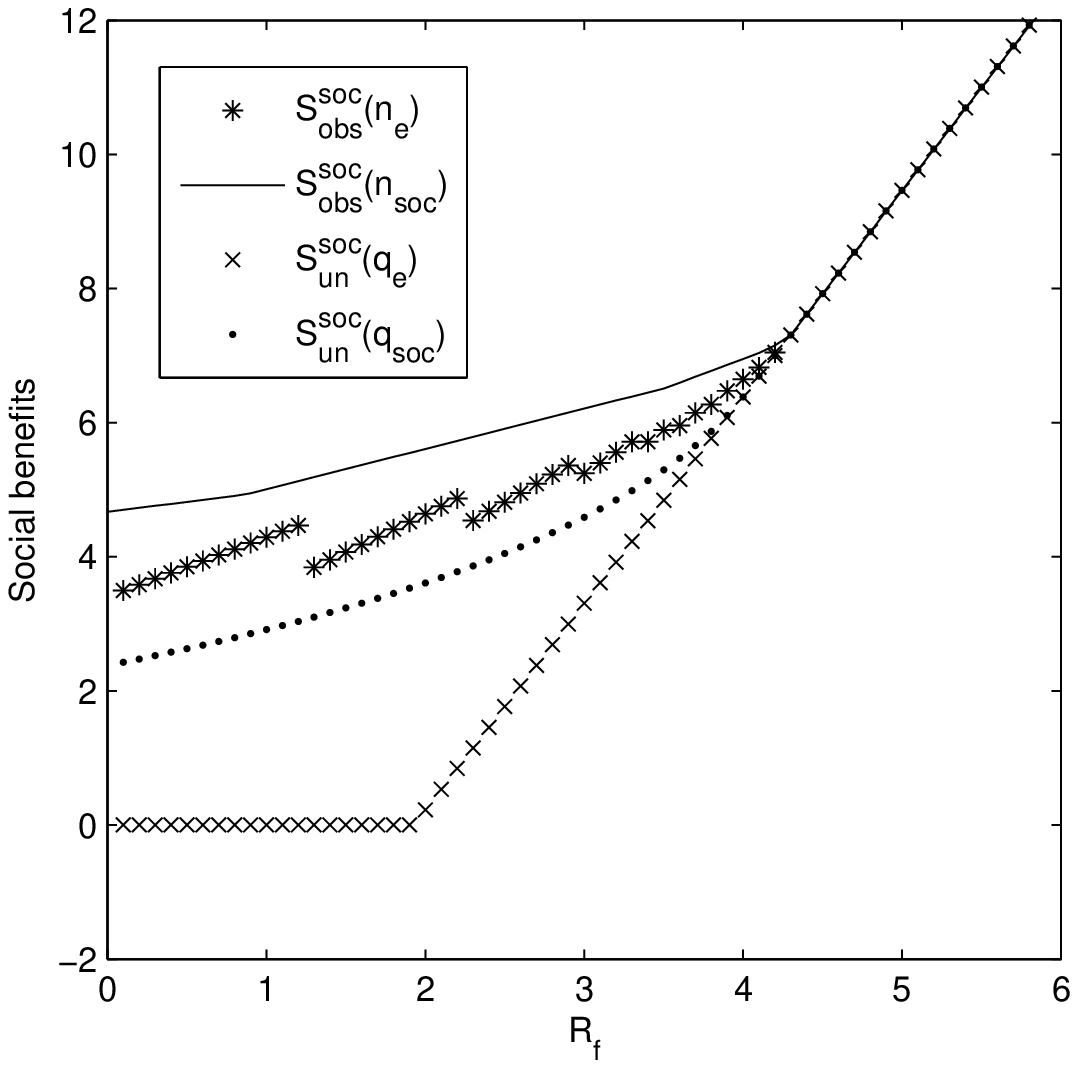}
\centering \caption{Social benefit per time unit with respect to
$R_f$ for a model with $(\lambda,\mu,\xi,\eta)=(7,2,0.7,1)$ and
$(R_s,C)=(7,3)$.} \label{numerical-scenario2}\end{figure}

\begin{figure}[hp]
\includegraphics[height=9cm,width=12cm,clip=true]{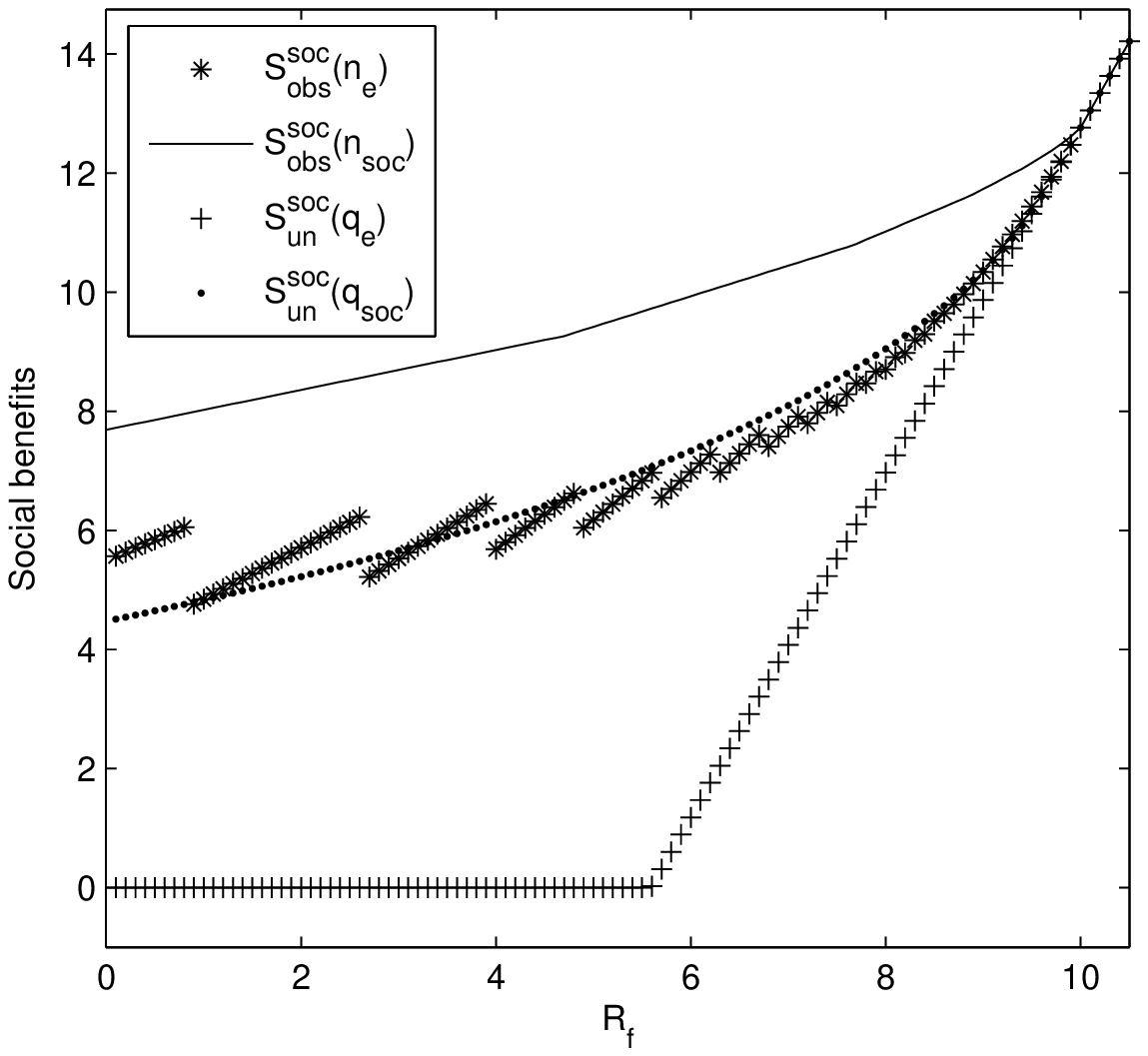}
\centering \caption{Social benefit per time unit with respect to
$R_f$ for a model with $(\lambda,\mu,\xi,\eta)=(7,4,0.3,2)$ and
$(R_s,C)=(4,3)$.} \label{numerical-scenario3}\end{figure}

Another topic of interest is the comparison between the observable
and the unobservable cases of a given model, i.e. what is the effect
of the information on customers' behavior. The value of the
information has been studied in a number of papers, among them in
Hassin (1986, 2007) and Guo and Zipkin (2007). In the context of the
present model, we have run several numerical scenarios and have
compared the expected social profit per time unit for the observable
and the unobservable cases, when the customers use their individual
or social optimal strategy. More concretely, we have been interested
in comparing $S_{obs}^{soc}(n_e)$, $S_{obs}^{soc}(n_{soc})$,
$S_{un}^{soc}(q_e)$ and $S_{un}^{soc}(q_{soc})$. The inequalities
$S_{obs}^{soc}(n_e) \leq S_{obs}^{soc}(n_{soc})$ and
$S_{un}^{soc}(q_e) \leq S_{un}^{soc}(q_{soc})$ are obviously valid
but the other relations are not clear. For example it would be
interesting to know the relationship between $S_{un}^{soc}(q_{soc})$
and $S_{obs}^{soc}(n_e)$ which corresponds to the natural question
`what is preferable for the society: to have uninformed altruistic
or informed selfish agents?'. The analysis of a large number of
numerical scenarios suggests that the typical ordering is
$S_{un}^{soc}(q_e) \leq S_{un}^{soc}(q_{soc}) \leq S_{obs}^{soc}
(n_e) \leq S_{obs}^{soc}(n_{soc})$. In this sense, it seems that in
the majority of such models it is better for the society the
customers to be informed and selfish than uninformed and altruistic.
Such a typical case is presented in Figure \ref{numerical-scenario2}
for $(\lambda,\mu,\xi,\eta)=(7,2,0.7,1)$ and $(R_s,C)=(7,3)$, as
$R_f$ varies in $[0,6]$. However, there are some exceptional cases
where for low values of $R_f$ we have that $S_{un}^{soc}(q_{soc})
\leq S_{obs}^{soc}(n_e)$, whereas for high values of $R_f$ we have
the reverse inequality. For intermediate values of $R_f$ the
situation is mixed. These cases occur typically for low values of $\xi$.
Such a numerical scenario is presented in Figure
\ref{numerical-scenario3} for $(\lambda,\mu,\xi,\eta)=(7,4,0.3,2)$
and $(R_s,C)=(4,3)$, as $R_f$ varies in $[0,10]$. Note also that
all graphs of $S_{obs}^{soc}(n_e)$, $S_{obs}^{soc}(n_{soc})$,
$S_{un}^{soc}(q_e)$ and $S_{un}^{soc}(q_{soc})$ with respect to $R_f$ 
coincide for $R_f\geq \frac{C}{\xi}$. Indeed, for values of $R_f$ 
exceeding the mean waiting cost till the next catastrophe, $\frac{C}{\xi}$,
it is both individually and socially optimal for the customers to 
join under any kind of information.

\section{Bibliography}

\begin{enumerate}

\item Burnetas, A. and Economou, A. (2007) Equilibrium customer
strategies in a single server Markovian queue with setup times.
\textit{Queueing Systems} \textbf{56}, 213-228.

\item Economou, A. and Kanta, S. (2008a) Optimal balking strategies
and pricing for the single server Markovian queue with compartmented
waiting space. \textit{Queueing Systems} \textbf{59}, 237-269.

\item Economou, A. and Kanta, S. (2008b) Equilibrium balking strategies
in the observable single-server queue with breakdowns and repairs.
\textit{Operations Research Letters} \textbf{36}, 696-699.

\item Edelson, N.M. and Hildebrand, K. (1975) Congestion tolls
for Poisson queueing processes. \textit{Econometrica} \textbf{43},
81-92.

\item Elaydi, S.N. (1999) \textit{An Introduction to Difference
Equations, 2nd Edition.} Springer, New York.

\item Guo, P. and Zipkin, P. (2007) Analysis and comparison of queues
with different levels of delay information. \textit{Management
Science} \textbf{53}, 962-970.

\item Hassin, R. (1986) Consumer information in markets with random
products quality: The case of queues and balking.
\textit{Econometrica} \textbf{54}, 1185-1195.

\item Hassin, R. (2007) Information and uncertainty in a queuing system.
\textit{Probability in the Engineering and Informational Sciences}
\textbf{21}, 361-380.

\item Hassin, R. and Haviv, M. (1995) Equilibrium strategies for
queues with impatient customers. \textit{Operations Research
Letters} \textbf{17}, 41-45.

\item Hassin, R. and Haviv, M. (1997) Equilibrium threshold strategies:
the case of queues with priorities. \textit{Operations Research}
\textbf{45}, 966-973.

\item Hassin, R. and Haviv, M. (2003) \textit{To Queue or Not to
Queue: Equilibrium Behavior in Queueing Systems.} Kluwer Academic
Publishers, Boston.

\item Mandelbaum, A. and Shimkin, N. (2000) A model for rational
abandonments from invisible queues. \textit{Queueing Systems}
\textbf{36}, 141-173.

\item Naor, P. (1969) The regulation of queue size by levying tolls.
\textit{Econometrica} \textbf{37}, 15-24.

\item Stidham, S. Jr. (2009) \textit{Optimal Design of Queueing
Systems.} CRC Press, Taylor and Francis Group, Boca Raton.

\end{enumerate}

\end{document}